\documentclass[12pt]{amsart}
\usepackage[bottom=3cm, right=3cm, left=3cm, top=3cm]{geometry}
\usepackage{amssymb,amsmath,amsthm,mathtools,amsfonts,times,hyperref,mathrsfs,multirow}
\usepackage{pgfplots}
\pgfplotsset{compat=1.15}
\usepackage{mathrsfs,mathdots}
\usetikzlibrary{arrows}
\usepackage{enumerate}
\usepackage{graphicx}
\usepackage{array}
\usepackage{ytableau}
\usepackage{pstricks,pst-node,graphicx}
\usepackage{tikz,varwidth}
\usepackage{setspace}
\usepackage{float}
\usetikzlibrary{calc}
\usepackage[super]{nth}
\usepackage{tikz}
\usepackage{extarrows}
\usepackage{pgfplots}
\pgfplotsset{compat=1.15}
\usetikzlibrary{arrows}

\newtheorem{theorem}{Theorem}[section]

\newtheorem{cor}[theorem]{Corollary}

\theoremstyle{definition}

\theoremstyle{remark}
\newtheorem{remark}{Remark}
\numberwithin{equation}{section}

\allowdisplaybreaks

\newcommand\numberthis{\addtocounter{equation}{1}\tag{\theequation}}
\newcommand\restr[2]{{
  \left.\kern-\nulldelimiterspace 
  #1 
  \littletaller 
  \right|_{#2} 
  }}
\newcommand{\littletaller}{\mathchoice{\vphantom{\big|}}{}{}{}}  

\begin{document}
\title[On Glaisher's Partition Theorem]{On Glaisher's Partition Theorem}
\author{George E. Andrews}
\address{Department of Mathematics, The Pennsylvania State University, University Park, PA 16802, USA}
\email{gea1@psu.edu}
\author{Aritram Dhar}
\address{Department of Mathematics, University of Florida, Gainesville, FL 32611, USA}
\email{aritramdhar@ufl.edu}

\date{\today}



\begin{abstract}
Glaisher's theorem states that the number of partitions of $n$ into parts which repeat at most $m-1$ times is equal to the number of partitions of $n$ into parts which are not divisible by $m$. The $m=2$ case is Euler's famous partition theorem. Recently, Andrews, Kumar, and Yee gave two new partition functions $C(n)$ and $D(n)$ related to Euler's theorem. Lin and Zang extended their result to Glaisher's theorem by generalizing $C(n)$. We generalize $D(n)$ and prove an analogous partition identity for the $m=3$ case. We also provide a new series equal to Glaisher's product both in the finite and infinite cases. 
\end{abstract}

\maketitle

\section{Introduction}\label{s1}
Euler's celebrated partition theorem \cite{A} can be stated as\\
\begin{align*}
A(n) = B(n)\numberthis\label{eq11}\\
\end{align*}
where $A(n)$ is the number of partitions of $n$ into distinct parts and $B(n)$ is the number of partitions of $n$ into odd parts.\\\quad\par In a recent work, Andrews, Kumar, and Yee \cite{AKY} have provided two more elegant companion partition functions to $A(n)$ and $B(n)$. In particular, they have obtained the following partition identity.\\
\begin{theorem}(\cite[Theorem $1.1$]{AKY})\label{thm11}
For any positive integer $n$, we have\\
\begin{align*}
A(n) = B(n) = C(n+1) = \dfrac{1}{2}D(n+1),\numberthis\label{eq12}\\
\end{align*}
where $C(n)$ is the number of partitions of $n$ with largest part even and parts not exceeding half of the largest part are distinct and $D(n)$ is the number of partitions of $n$ into non-negative parts wherein the smallest part appear exactly twice and no other parts are repeated.\\
\end{theorem}
Let $m\ge 2$ be an integer. Define $A_m(n)$ to be the number of partitions of n into parts repeating less than $m$ times and $B_m(n)$ to be the number of partitions of $n$ into parts not divisible by $m$. In 1883, Glaisher \cite{G} obtained the following partition theorem which is a natural extension of Euler's theorem \eqref{eq11}.\\  
\begin{theorem}(\cite{G})\label{thm12}
For any non-negative integer $n$, we have\\
\begin{align*}
A_m(n) = B_m(n).\numberthis\label{eq13}\\
\end{align*}
\end{theorem}
For any integer $0\le j\le m-1$, define $B^{(j)}_m(n)$ to be the number of partitions of $n$ into parts not divisible by $m$ and largest part congruent to $j$ modulo $m$. Clearly,\\
\begin{align*}
B_m(n) = \sum\limits_{j=1}^{m-1}B^{(j)}_m(n).\numberthis\label{eq14}\\
\end{align*}
\par In a recent work, Lin and Zang \cite{LZ} generalized Andrews, Kumar, and Yee's $C(n)$ function to the case of $m$-\textit{regular partitions} (partitions whose parts are not divisible by $m$). It is noteworthy to point out that the cardinality of the set $B^{\prime}_m(n)$ in their paper is the same as $B^{(m-1)}_m(n)$.\\\quad\par Define $C_m(n)$ to be the number of partitions of $n$ in which the largest part is divisible by $m$ (say $mj$) and all parts less than or equal to $j$ repeat less than $m$ times. Also, define $C_m(0) := 1$.\\
\begin{theorem}(\cite[Theorem $1.5$]{LZ})\label{thm13}
For any positive integer $n$, we have\\
\begin{align*}
B^{(m-1)}_m(n) = C_m(n+1).\numberthis\label{eq15}\\
\end{align*}
\end{theorem}
\begin{remark}\label{rmk1}
Lin and Zang \cite{LZ} provided both $q$-series (through generating functions) and bijective proofs of Theorem \ref{thm13}.\\
\end{remark}
\par Let us now recall some notations from the theory of $q$-series which can be found in \cite{A}. For any complex variable $A$ and non-negative integer $N$, the conventional $q$-Pochhammer symbol is defined as\\
\begin{align*}
(A)_N = (A;q)_N &:= \Bigg\{\begin{array}{lr}
1\qquad\qquad\quad\,\,\,\text{when } N=0,\\
\prod\limits_{i=0}^{N-1}(1-Aq^i)\,\,\,\text{otherwise, and}\end{array}\\
(A)_{\infty} = (A;q)_{\infty} &:= \lim\limits_{N\rightarrow\infty}(A)_N\,\,\text{for}\,\,|q|<1.\\
\end{align*}
Also, let $\zeta_n = e^{\frac{2\pi i}{n}}$ denote the primitive $n^{\text{th}}$ root of unity and $\omega = \zeta_3$ denote the primitive cubic root of unity here and throughout the remainder of the paper.\\\quad\par Now, define $D_m(n)$ to be the number of partitions of $n$ into non-negative parts where the smallest part occurs exactly $m$ times and all other parts repeat less than $m$ times.\\
\begin{theorem}\label{thm14}
For any non-negative integer $n$, we have\\
\begin{align*}
C_m(n) = \dfrac{D_m(n) + E_m(n)}{m},\numberthis\label{eq16}\\
\end{align*}
where $E_m(n)$ is such that the generating function of $E_m(n)$ is\\
\begin{align*}
\varepsilon_m(q) := \sum\limits_{n=0}^{\infty}E_m(n)q^n = \sum\limits_{n=0}^{\infty}q^{mn}(q^{m(n+1)};q)_{\infty}\sum\limits_{j=1}^{m-1}\dfrac{1}{(\zeta^{j}_mq^{n+1};q)_{\infty}}.\numberthis\label{eq17}\\
\end{align*}
\end{theorem}
The $m = 3$ case of \eqref{eq17} then corresponds to the following identity.\\
\begin{theorem}\label{thm15}
\begin{align*}
\varepsilon_3(q) := \sum\limits_{n=0}^{\infty}E_3(n)q^n = 2-q-2q^2+\sum\limits_{n=2}^{\infty}(-1)^n\chi(n-1)q^{\binom{n+1}{2}+1}.\numberthis\label{eq18}\\    
\end{align*}
\end{theorem}
The complexity of the generating functions involved seemed to mitigate against an analog of the $C(n) = \frac{1}{2}D(n)$ portion of Theorem \ref{thm11}; however, surprisingly, a similar result exists for $m = 3$ which clearly follows from \eqref{eq16} and \eqref{eq18}.\\
\begin{cor}\label{cor16}
Let $n$ be any positive integer that is not equal to a triangular number plus $1$, then\\
\begin{align*}
C_3(n) = \frac{1}{3}D_3(n).\numberthis\label{eq19}\\ 
\end{align*}
\end{cor}
Combining Theorems \ref{thm12}, \ref{thm13}, and \ref{thm14}, we then have the following partition identity which is a natural extension of Theorem \ref{thm11}.\\
\begin{theorem}\label{thm17}
For any positive integer $n$, we have\\
\begin{align*}
A_m(n) = B_m(n) = \sum\limits_{k=1}^{m-2}B^{(k)}_m(n) + C_m(n+1) = \sum\limits_{k=1}^{m-2}B^{(k)}_m(n) + \dfrac{D_m(n+1) + E_m(n+1)}{m}.\numberthis\label{eq110}\\
\end{align*}
\end{theorem}
In view of the generating function of $B_2(n)$ in \cite[$(2.1)$]{AKY}, we also obtain the following finite $q$-series identity.\\
\begin{theorem}\label{thm18}
For any positive integer $N$, we have\\
\begin{align*}
1+\sum\limits_{n=1}^{N}\sum\limits_{j=1}^{m-1}\dfrac{q^{mn-j}}{\left(\prod\limits_{r=1}^{m-j}(q^r;q^m)_n\right)\cdot\left(\prod\limits_{r=m-j+1}^{m-1}(q^r;q^m)_{n-1}\right)} = \dfrac{(q^m;q^m)_{N}}{(q;q)_{mN}}.\numberthis\label{eq111}\\
\end{align*}
\end{theorem}
Letting $N\rightarrow\infty$ in \eqref{eq111} gives us the following $q$-series identity.\\
\begin{cor}\label{cor19}
\begin{align*}
1+\sum\limits_{n=1}^{\infty}\sum\limits_{j=1}^{m-1}\dfrac{q^{mn-j}}{\left(\prod\limits_{r=1}^{m-j}(q^r;q^m)_n\right)\cdot\left(\prod\limits_{r=m-j+1}^{m-1}(q^r;q^m)_{n-1}\right)} = \dfrac{(q^m;q^m)_{\infty}}{(q;q)_{\infty}}.\numberthis\label{eq112}\\    
\end{align*}    
\end{cor}
\begin{remark}\label{rmk2}
The $m=2$ case of \eqref{eq112} above was recently proved by Andrews and Bachraoui in \cite[Lemma $1$, (1.10)]{AB}. Clearly, the series on the left-hand side of \eqref{eq112} is the generating function of $\sum\limits_{j=1}^{m-1}B^{(j)}_m(n)$ and the product on the right-hand side of \eqref{eq112} is the generating function of $B_m(n)$. The same observation (restricting all parts to be at most $mN$) proves Theorem \ref{thm18}.\\ 
\end{remark}
The rest of the paper is organized as follows. In Sections \ref{s2} and \ref{s3}, we provide $q$-series (generating function) proofs of Theorem \ref{thm14} and Theorem \ref{thm15} respectively.\\

\section{Proof of Theorem \ref{thm14}}\label{s2}
From the definition of $C_m(n)$, we then have\\
\begin{align*}
\sum\limits_{n=0}^{\infty}C_m(n)q^n &= \sum\limits_{n=0}^{\infty}\dfrac{\prod\limits_{j=1}^{n}(1+q^{j}+\ldots+q^{(m-1)j})q^{mn}}{(q^{n+1};q)_{(m-1)n}}\\
&= \sum\limits_{n=0}^{\infty}\dfrac{(q^m;q^m)_nq^{mn}}{(q;q)_{mn}}\\
&= (q^m;q^m)_{\infty}\sum\limits_{n=0}^{\infty}\dfrac{q^{mn}}{(q;q)_{mn}(q^{mn+m};q^m)_{\infty}}.\numberthis\label{eq21}\\
\end{align*}
Using Euler's identity \cite[p. $19$, ($2.2.5$)]{A}\\
\begin{align*}
\dfrac{1}{(t;q)_{\infty}} = \sum\limits_{n=0}^{\infty}\dfrac{t^n}{(q;q)_n},\numberthis\label{eq22}\\    
\end{align*}
(replacing $q$ by $q^m$ and letting $t=q^{mn+m}$) in \eqref{eq21}, we get\\
\begin{align*}
\sum\limits_{n=0}^{\infty}C_m(n)q^n &= (q^m;q^m)_{\infty}\sum\limits_{n=0}^{\infty}\dfrac{q^{mn}}{(q;q)_{mn}}\sum\limits_{j=0}^{\infty}\dfrac{q^{mnj+mj}}{(q^m;q^m)_j}\\
&= (q^m;q^m)_{\infty}\sum\limits_{n,j=0}^{\infty}\dfrac{q^{mn+mnj+mj}}{(q;q)_{mn}(q^m;q^m)_j}\\
&= (q^m;q^m)_{\infty}\sum\limits_{n,j=0}^{\infty}\dfrac{1}{m}\left(\sum\limits_{k=0}^{m-1}\zeta^{kn}_m\right)\dfrac{q^{n+nj+mj}}{(q;q)_{n}(q^m;q^m)_j}\\
&= \dfrac{(q^m;q^m)_{\infty}}{m}\sum\limits_{j=0}^{\infty}\dfrac{q^{mj}}{(q^m;q^m)_j}\sum\limits_{k=0}^{m-1}\sum\limits_{n=0}^{\infty}\dfrac{(\zeta^k_mq^{j+1})^n}{(q;q)_n}\numberthis\label{eq23}.\\
\end{align*}
Using \eqref{eq22} with $t=\zeta^k_mq^{j+1}$ in \eqref{eq23}, we get\\
\begin{align*}
\sum\limits_{n=0}^{\infty}C_m(n)q^n &= \dfrac{(q^m;q^m)_{\infty}}{m}\sum\limits_{j=0}^{\infty}\dfrac{q^{mj}}{(q^m;q^m)_j}\sum\limits_{k=0}^{m-1}\dfrac{1}{(\zeta^k_mq^{j+1};q)_{\infty}}\\
&= \dfrac{(q^m;q^m)_{\infty}}{m}\sum\limits_{j=0}^{\infty}\dfrac{q^{mj}}{(q^m;q^m)_j(q^{j+1};q)_{\infty}}\\
&\qquad\qquad + \dfrac{(q^m;q^m)_{\infty}}{m}\sum\limits_{j=0}^{\infty}\dfrac{q^{mj}}{(q^m;q^m)_j}\sum\limits_{k=1}^{m-1}\dfrac{1}{(\zeta^k_mq^{j+1};q)_{\infty}}\\
&= \dfrac{1}{m}\sum\limits_{j=0}^{\infty}q^{mj}\cdot\prod\limits_{i=j+1}^{\infty}(1+q^i+\ldots+q^{(m-1)i})\\
&\qquad\qquad + \dfrac{1}{m}\sum\limits_{j=0}^{\infty}q^{mj}(q^{m(j+1)};q)_{\infty}\sum\limits_{k=1}^{m-1}\dfrac{1}{(\zeta^k_mq^{j+1};q)_{\infty}}\\
&= \dfrac{1}{m}\sum\limits_{n=0}^{\infty}D_m(n)q^n + \dfrac{1}{m}\sum\limits_{n=0}^{\infty}E_m(n)q^n,\numberthis\label{eq24}\\
\end{align*}
where the generating function of $D_m(n)$ clearly follows from its definition and the generating function of $E_m(n)$ follows from \eqref{eq17}. Comparing coefficients of $q^n$ on both sides of \eqref{eq24} for any non-negative integer $n$, we have\\
\begin{align*}
C_m(n) = \dfrac{D_m(n) + E_m(n)}{m}\\   
\end{align*}
which completes the proof of Theorem \ref{thm14}.\\\qed\\

\section{Proof of Theorem \ref{thm15}}\label{s3}
For $m=3$, \eqref{eq17} can be written as\\
\begin{align*}
\varepsilon_3(q) = (q;q)_{\infty}\sum\limits_{n=0}^{\infty}\dfrac{q^{3n}}{(q;q)_n}\left((\omega q^{n+1};q)_{\infty}+(\omega^{-1}q^{n+1};q)_{\infty}\right).\numberthis\label{eq31}   
\end{align*}
Define $\chi(n)$ as follows\\
\begin{align*}
\chi(n) := \omega^n + \omega^{-n} = \Bigg\{\begin{array}{lr}
\,\,\,\,2\,\,\,\,\,\text{if } 3\mid n,\\
-1\,\,\,\,\,\text{if } 3\nmid n.\end{array}\\
\end{align*}
Using \cite[p. $19$, ($2.2.6$)]{A} with $t = -\omega^jq^n$ in \eqref{eq31}, we get\\
\begin{align*}
\varepsilon_3(q) &= (q;q)_{\infty}\sum\limits_{n=0}^{\infty}\dfrac{q^{3n}}{(q;q)_n}\sum\limits_{k=0}^{\infty}\dfrac{(-1)^k\chi(k)q^{\binom{k+1}{2}+nk}}{(q;q)_k}.\numberthis\label{eq32}\\
\end{align*}
Using Euler's identity \cite[p. $19$, ($2.2.5$)]{A} with $t = q^{k+3}$ after interchanging summations in \eqref{eq32}, we then have\\
\begin{align*}
\varepsilon_3(q) &= (q;q)_{\infty}\sum\limits_{k=0}^{\infty}\dfrac{(-1)^k\chi(k)q^{\binom{k+1}{2}}}{(q;q)_k(q^{k+3};q)_{\infty}}\\
&= \sum\limits_{k=0}^{\infty}(-1)^k\chi(k)q^{\binom{k+1}{2}}(1-q^{k+1})(1-q^{k+2})\\
&= \sum\limits_{k=0}^{\infty}(-1)^k\chi(k)q^{\binom{k+1}{2}}(1-q^{k+1}(1+q)+q^{(k+1)+(k+2)})\\
&= \sum\limits_{k=0}^{\infty}(-1)^k\chi(k)q^{\binom{k+1}{2}} - (1+q)\sum\limits_{k=0}^{\infty}(-1)^k\chi(k)q^{\binom{k+2}{2}} + \sum\limits_{k=0}^{\infty}(-1)^k\chi(k)q^{\binom{k+3}{2}}\\
&= \sum\limits_{k=0}^{\infty}(-1)^k\chi(k)q^{\binom{k+1}{2}} + (1+q)\sum\limits_{k=1}^{\infty}(-1)^k\chi(k-1)q^{\binom{k+1}{2}} + \sum\limits_{k=2}^{\infty}(-1)^k\chi(k-2)q^{\binom{k+1}{2}}\\
&= 2 + q + \sum\limits_{k=2}^{\infty}(-1)^k\chi(k)q^{\binom{k+1}{2}} - 2q(1+q) + (1+q)\sum\limits_{k=2}^{\infty}(-1)^k\chi(k-1)q^{\binom{k+1}{2}}\\
&\qquad\quad + \sum\limits_{k=2}^{\infty}(-1)^k\chi(k-2)q^{\binom{k+1}{2}}\\
&= 2 - q - 2q^2 + \sum\limits_{k=2}^{\infty}(-1)^k\left(\chi(k)+\chi(k-1)+\chi(k-2)\right)q^{\binom{k+1}{2}}\\
&\qquad\quad + \sum\limits_{k=2}^{\infty}(-1)^k\chi(k-1)q^{\binom{k+1}{2}+1}.\\
\end{align*}
Note that $\chi(k)+\chi(k-1)+\chi(k-2) = 0$ for all $k$. Hence, we have\\
\begin{align*}
\varepsilon_3(q) = 2 - q - 2q^2 + \sum\limits_{k=2}^{\infty}(-1)^k\chi(k-1)q^{\binom{k+1}{2}+1}.\\     
\end{align*}
This completes the proof of Theorem \ref{thm15}.\\\qed\\

\section{Concluding Remarks}\label{s4}
It would be very interesting to provide a bijective proof of Corollary \ref{cor16}.\\

\section*{Acknowledgments}\label{s5}
The authors would like to thank Alexander Berkovich and Rahul Kumar for their kind interest. The authors would also like to thank the anonymous referee for their helpful comments and suggestions.\\

\end{document}